\documentclass[11pt]{article}
\usepackage[tbtags]{amsmath}
\usepackage{amsfonts,verbatim,amssymb,amsthm}
\usepackage{mathrsfs}
\usepackage{graphicx}
\allowdisplaybreaks

\numberwithin{equation}{section}

\theoremstyle{plain}
\newtheorem{theorem}{Theorem}

\newtheorem{corollary}{Corollary}
\newtheorem{proposition}{Proposition}

\theoremstyle{definition}

\newtheorem{remark}{Remark}

\begin{document}


\title{On construction of transmutation operators for perturbed
Bessel equations}

\author{Vladislav V. Kravchenko$^1$, Elina L. Shishkina$^2$ and Sergii M. Torba$^1$\\
{\small $^1$ Departamento de Matem\'{a}ticas, CINVESTAV del IPN, Unidad
Quer\'{e}taro,}\\
{\small Libramiento Norponiente \#2000, Fracc.\ Real de Juriquilla, Quer\'{e}taro, Qro., 76230 MEXICO.} \\
{\small $^2$ Faculty of Applied
Mathematics, Informatics and Mechanics,}\\
{\small Voronezh State University,
Universitetskaya pl., 1, 394080, Voronezh, RUSSIA.}\\
{\small \texttt{vkravchenko@math.cinvestav.edu.mx},}\\
{\small\texttt{shishkina@amm.vsu.ru},}\\
{\small\texttt{storba@math.cinvestav.edu.mx}}
}

\date{December~04, 2017}


\markboth{V.~V.~Kravchenko, E.~L.~Shishkina, S.~M.~Torba}{On construction of
transmutation operators for perturbed Bessel equations}

\maketitle

\begin{abstract}
A representation for the kernel of the transmutation operator relating the
perturbed Bessel equation with the unperturbed one is obtained in the form of a
functional series with coefficients calculated by a recurrent integration
procedure. New properties of the transmutation kernel are established. A new
representation of the regular solution of the perturbed Bessel equation is
given presenting a remarkable feature of uniform error bound with respect to
the spectral parameter for partial sums of the series. A numerical illustration
of application to the solution of Dirichlet spectral problems is presented.
\end{abstract}


\section{Introduction}\label{sec:level1}

In \cite{Volk} it was proved that a regular solution of the perturbed Bessel
equation
\begin{equation}
-u^{\prime \prime }+\left( \frac{l(l+1)}{x^{2}}+q(x)\right) u=\omega ^{2}u,\quad x\in (0,b],  \label{Eq01}
\end{equation}%
where $q\in C[0,b]$, $l\geq -\frac{1}{2}$, for all $\omega \in \mathbb{C}$ can
be obtained from a regular solution of the unperturbed Bessel equation
\begin{equation}
-y^{\prime \prime }+\frac{l(l+1)}{x^{2}}y=\omega ^{2}y  \label{BesselEq}
\end{equation}%
with the aid of a transmutation (transformation) operator in the form of a
Volterra integral operator,
\begin{equation}
u(\omega ,x)=T[y(\omega ,x)]=y(\omega ,x)+\int\limits_{0}^{x}K(x,t)y(\omega ,t)dt.  \label{transmutation operator}
\end{equation}%
Here the kernel $K$ is $\omega $-independent continuous function with respect
to both arguments satisfying the Goursat condition
\begin{equation}\label{Goursat}
    K(x,x) = \frac 12\int\limits_0^x q(s)\,ds.
\end{equation}
A regular solution of (\ref{BesselEq}) can be chosen in the form
\begin{equation*}
y(\omega ,x)=\sqrt{x}J_{l+\frac{1}{2}}(\omega x)
\end{equation*}%
where $J_{\nu }$ stands for the Bessel function of the first kind and of order
$\nu $.
The transmutation operator (\ref{transmutation operator}) is a fundamental
object in the theory of inverse problems related to (\ref{Eq01}) and has been
studied in a number of publications (see, e.g., \cite{Chadan}, \cite{Chebli},
\cite{Coz1976}, \cite{Kostenko_dispersionestimates}, \cite{Sta2},
\cite{Sitnik},  \cite{Trimeche}).


Up to now apart from a successive approximation procedure used in \cite{Volk},
\cite{Coz1976} for proving the existence of $K$ and a series representation
proposed in \cite{Chebli} requiring the potential to possess holomorphic
extension onto the disk of radius $2xe\sqrt{1+|2l|}$, no other construction of
the transmutation kernel $K$ has been proposed. In this relation we mention the
recent work \cite{KTS} where $K$ was approximated by a special system of
functions called generalized wave polynomials.

In the present paper we obtain an exact representation of the kernel $K$ in the
form of a functional series whose coefficients are calculated following a
simple recurrent integration procedure. The representation has an especially
simple and attractive form in the case when $l$ is a natural number. It
revealed some new properties of the kernel $K$.

The representation is obtained with the aid of a recent result from \cite%
{KTC2017} where a Fourier-Legendre series expansion was derived for a certain
kernel $R$ related to the kernel $K$. Here with the aid of an Erdelyi-Kober
fractional derivative we express $K$ in terms of $R$ which leads to a series
representation for the kernel $K$.

The obtained form of the kernel $K$ is appropriate both for studying the exact
solution as well as the properties of the kernel itself, and for numerical
applications.

In Section \ref{sec:level2} some previous results are recalled and an
expression for the kernel $K$ in terms of $R$ is derived. In Section
\ref{Sect3} a functional series representation for $K$ is obtained in the case
of integer parameter $l$. Its convergence properties are studied.

In Section  \ref{Sect4} the obtained representation of the kernel $K$ is used
for deriving a new representation for regular solution of the perturbed Bessel
equation enjoying the uniform ($\omega$-independent) approximation property
(Theorem \ref{Thm uniform}). Since the regular solution $u(\omega, x)$ in this
representation is the image of $\sqrt{\omega x} J_{l+1/2}(\omega x)$ under the
action of the transmutation operator $T$, $u(\omega, x)$ does not decay to zero
as $\omega\to \infty$. Hence partial sums of the representation provide good
approximation even for arbitrarily large values of the spectral parameter
$\omega$. This is an advantage in comparison with the representation from
\cite{KTC2017} derived for a regular solution decaying as $\omega^{-l-1}$ when
$\omega\to\infty$, making the approximation by partial sums useful for
reasonably small values of $\omega$ only.

In Section \ref{Sect5} a functional series representation for the kernel $K$ is
obtained in the case of a noninteger $l$. Section \ref{Sect6} contains some
numerical illustrations confirming the validity of the presented results.

\section{A representation of $K(x,t)$ using an Erdelyi-Kober
operator}\label{sec:level2}

In this section, with the aid of a result from \cite{KTC2017} we obtain a
representation of the kernel $K(x,t)$ in terms of an Erdelyi-Kober fractional
derivative applied to a certain Fourier-Legendre series. Two transmutation
operators are used, a modified Poisson transmutation operator and the
transmutation operator (\ref{transmutation operator}).

\subsection{The modified Poisson transmutation operator}\label{Subsect The
modified Poisson}

The Poisson transmutation operator examined in \cite{levitan} and \cite%
{Sitnik} is adapted to work with the singular differential Bessel operator $%
B_{\gamma }=\frac{d^{2}}{dx^{2}}+\frac{\gamma }{x}\frac{d}{dx}$, where $%
\gamma >0$. We will use a slightly modified Poisson operator defined on $%
C([0,b])$ of the form (see \cite{KTC2017} and \cite{KTS})
\begin{equation}
Y_{l,x}f(x)=\frac{x^{-l}}{2^{l+\frac{1}{2}}\Gamma \left( l+\frac{3}{2}%
\right) }\int\limits_{0}^{x}(x^{2}-s^{2})^{l}f(s)ds,\quad l\geq -\frac{1}{2}.
 \label{Poisson}
\end{equation}%
The following equality is valid
\begin{equation*}
Y_{l,x}[\cos {\omega x}]=\frac{\sqrt{\pi }\Gamma (l+1)}{2\omega ^{l+1}\Gamma \left( l+\frac{3}{2}\right) }\,\sqrt{\omega x}\,J_{l+\frac{1}{2}}(\omega x).
\end{equation*}%
Moreover, $Y_{l,x}$ is an intertwining operator for $\frac{d^{2}}{dx^{2}}$
and $\frac{d^{2}}{dx^{2}}-\frac{l(l+1)}{x^{2}}$ in the following sense. If $%
v\in C^{2}[0,b]$ and $v^{\prime }(0)=0$ then
\begin{equation*}
Y_{l,x}\frac{d^{2}v}{dx^{2}}=\left( \frac{d^{2}}{dx^{2}}-\frac{l(l+1)}{x^{2}}%
\right) Y_{l,x}v.
\end{equation*}%
In particular, the regular solution of equation \eqref{BesselEq} satisfying
the asymptotic relations $y_{0}(x){\sim }x^{l+1}$ and $y_{0}^{\prime }(x){%
\sim }(l{+}1)x^{l}$, when $x\rightarrow 0$ can be written in the form
\begin{equation}
y(\omega ,x)=\frac{2^{l+\frac{3}{2}}\Gamma ^{2}\left( l+\frac{3}{2}\right) }{%
\sqrt{\pi }\Gamma (l+1)}Y_{l,x}[\cos (\omega x)]=\Gamma \left( l+\frac{3}{2}%
\right) 2^{l+\frac{1}{2}}\omega ^{-l-\frac{1}{2}}\sqrt{x}J_{l+\frac{1}{2}%
}(\omega x).  \label{igrek}
\end{equation}

The composition of two transmutation operators (\ref{transmutation operator}%
) and \eqref{Poisson} allows one to write a regular solution of \eqref{Eq01}
in the form (see \cite{KTC2017})%
\begin{equation}
u(\omega ,x)=\frac{2^{l+\frac{1}{2}}\Gamma \left( l+\frac{3}{2}\right) }{%
\omega ^{l+\frac{1}{2}}}\sqrt{x}\,J_{l+\frac{1}{2}}(\omega x)
+\int\limits_{0}^{x}R(x,t)\cos {\omega t}dt,  \label{Sol}
\end{equation}
where the kernel $R(x,t)$ is a sufficiently regular function which admits a
convergent Fourier--Le\-gen\-dre series expansion presented in \cite{KTC2017}.
Namely,
\begin{equation}
R(x,t)=\sum\limits_{k=0}^{\infty }\frac{\beta _{k}(x)}{x}P_{2k}\left( \frac{t%
}{x}\right) ,  \label{Series1}
\end{equation}%
where $P_{m}(x)$ stands for the Legendre polynomial of order $m$, and the
coefficients $\beta _{k}$ can be computed following a simple recurrent
integration procedure \cite{KTC2017}.

\subsection{A representation of $K(x,t)$ in terms of $R(x,t)$}

\begin{theorem}
Let $q\in C[0,b]$ be a complex-valued function. Then the following equality is
valid
\begin{equation}
K(x,t)=\frac{\sqrt{\pi }}{\Gamma \left( l+\frac{3}{2}\right) }\frac{t^{l+1}}{\Gamma (n-l-1)} \left( -\frac{d}{2tdt}\right) ^{n}\int\limits_{t}^{x}(s^{2}-t^{2})^{n-l-2}sR(x,s)ds,  \label{Int4}
\end{equation}
here $n$ can be arbitrary integer satisfying $n>l+1$.
\end{theorem}

\noindent\textbf{Proof.} The following relation between the kernels $K$ and $R$
was obtained in \cite[(3.10) and (3.12)]{KTC2017}
\begin{equation}
\frac{2\Gamma \left( l+\frac{3}{2}\right) }{\sqrt{\pi }\Gamma (l+1)}%
\int\limits_{s}^{x}K(x,t)t^{-l}(t^{2}-s^{2})^{l}dt=R(x,s).  \label{R0}
\end{equation}%
Let us invert this equality using an inverse Erdelyi-Kober operator. For
$\alpha >0$ the left-sided Erdelyi-Kober operator is (see \cite[(18.3)]{SKM})
\begin{equation}
I_{x-;2,\eta }^{\alpha }f(s){=}\frac{2s^{\sigma \eta }}{\Gamma (\alpha )}%
\int\limits_{s}^{x}(t^{2}-s^{2})^{\alpha -1}t^{2(1-\alpha -\eta )-1}f(t)dt, \label{EK}
\end{equation}%
that means we can rewrite \eqref{R0} as
\begin{equation}
\frac{\Gamma \left( l+\frac{3}{2}\right) }{\sqrt{\pi }}\,t^{l+1}\,I_{x-;2,-%
\frac{l+1}{2}}^{l+1}K(x,t)=R(x,t).
\end{equation}%
Now applying the inverse operator for \eqref{EK} we obtain (\ref{Int4}).
\hfill\qed

\section{Representation of $K(x,t)$ for $l=0,1,2,...$}\label{Sect3}

The relation (\ref{Int4}) together with the Fourier-Legendre series
representations of $R$ in (\ref{Series1}) lead to a Fourier-Jacobi series
representation of the kernel $K$. It admits an especially simple form in the
case of integer values of the parameter $l$.


\begin{theorem}\label{Thm2}
Let $q\in C^1[0,b]$. For $l=0,1,2,...$ the following formula for the kernel of
the transmutation operator \eqref{transmutation operator} is valid
\begin{equation}
K(x,t)=\frac{\sqrt{\pi }t^{l+1}}{x^{2l+3}\Gamma \left( l+\frac{3}{2}\right) }
\sum\limits_{m=0}^{\infty }\frac{(-1)^{m+l+1}\Gamma \left( m+2l+\frac{5}{2}%
\right) }{\Gamma \left( m+l+\frac{3}{2}\right) } \beta _{m+l+1}(x)P_{m}^{\left( l+\frac{1}{2},l+1\right) }\left( 1-2\frac{t^{2}}{x^{2}}\right), \label{K(x,t)}
\end{equation}
where $\,P_{m}^{\left( \alpha ,\beta \right) }$ stands for a Jacobi
polynomial and the coefficients $\beta _{k}$ are those from Subsection \ref%
{Subsect The modified Poisson}.

The series \eqref{K(x,t)} converges absolutely for any $x\in (0,b]$ and $t\in
(0,x)$ and converges uniformly with respect to $t$ on any segment
$[\varepsilon, x-\varepsilon]\subset (0,x)$. Under the additional assumption
that $q\in C^{2l+5}[0,b]$ the series \eqref{K(x,t)} converges absolutely and
uniformly with respect to $t$ on $[0,x]$.
\end{theorem}

\noindent\textbf{Proof.}
For any nonnegative integer $l$ and $n>l+1$ we can rewrite formula %
\eqref{Int4} in the form
\begin{equation}
K(x,t)=\frac{(-1)^{l+1}\sqrt{\pi }}{\Gamma \left( l+\frac{3}{2}\right) }%
\frac{t^{l+1}}{\Gamma (n-l-1)} \left( \frac{d}{2tdt}\right) ^{n}\int\limits_{x}^{t}(t^{2}-s^{2})^{n-l-2}sR(x,s)ds.  \label{Int7}
\end{equation}
Choosing $n=l+2$ we obtain
\begin{equation}
K(x,t)=\frac{(-1)^{l+1}\sqrt{\pi }t^{l+1}}{\Gamma \left( l+\frac{3}{2} \right) }\left( \frac{d}{2tdt}\right) ^{l+2}\int\limits_{x}^{t}sR(x,s)ds=
\frac{(-1)^{l+1}\sqrt{\pi }t^{l+1}}{\Gamma \left( l+\frac{3}{2}\right) }%
\left( \frac{d}{2tdt}\right) ^{l+1}R(x,t).\\\label{Int71}
\end{equation}
Substituting the series expansion \eqref{Series1} into \eqref{Int71} we get (we
left the justification of the possibility of termwise differentiation of
\eqref{Series1} to the end of the proof)
\begin{equation*}
K(x,t)=\frac{(-1)^{l+1}\sqrt{\pi }t^{l+1}}{\Gamma \left( l+\frac{3}{2}%
\right) }\sum\limits_{k=0}^{\infty }\frac{\beta _{k}(x)}{x}\left( \frac{d}{%
2tdt}\right) ^{l+1}P_{2k}\left( \frac{t}{x}\right) .
\end{equation*}%
Now we would like to calculate the derivatives
\begin{equation*}
\left( \frac{d}{2tdt}\right) ^{l+1}P_{2k}\left( \frac{t}{x}\right) =\left(
\frac{d}{dt^{2}}\right) ^{l+1}P_{2k}\left( \frac{t}{x}\right).
\end{equation*}
 Note that by formula 8.911
from \cite{Gradshteyn} the following equality is valid
\begin{equation*}
\begin{split}
P_{2k}\left( \frac{t}{x}\right) &=
(-1)^{k}\frac{(2k-1)!!}{2^{k}k!}%
\,_{2}F_{1}\left( -k,k+\frac{1}{2};\frac{1}{2};\frac{t^{2}}{x^{2}}\right) 
= (-1)^{k}\frac{\left(\frac 12\right)_k}{k!} \,_{2}F_{1}\left(
-k,k+\frac{1}{2};\frac{1}{2};\frac{t^{2}}{x^{2}}\right) .
\end{split}
\end{equation*}%
Now application of formula 15.2.2 from \cite{AbramowitzStegunSpF},
\begin{equation*}
\frac{d^{n}}{dz^{n}}F(a,b;c;z)=\frac{(a)_{n}(b)_{n}}{(c)_{n}}F(a+n,b+n;c+n;z)
\end{equation*}%
leads to the relation
\begin{equation*}
\left( \frac{d}{dt^2}\right) ^{l+1}P_{2k}\left( \frac{t}{x}\right) =(-1)^{k}%
\frac{\left(\frac 12\right)_k}{k!}
\frac{(-k)_{l+1}\left( k+\frac{1}{2}\right) _{l+1}}{%
x^{2l+2}\left( \frac{1}{2}\right) _{l+1}}
\,_{2}F_{1}\left( l+1-k,k+l+\frac{3%
}{2};l+\frac{3}{2};\frac{t^{2}}{x^{2}}\right) .
\end{equation*}%
Taking into account that $(-k)_{l+1} = 0$ when $l\geq k$ and $(-k)_{l+1} =
(-1)^{l+1}\frac{k!}{\Gamma (k-l)}$ when $l<k$, we obtain
\begin{equation*}
\left( \frac{d}{dt^2}\right) ^{l+1}P_{2k}\left( \frac{t}{x}\right) =0
\end{equation*}%
when $l\geq k$, and
\begin{equation*}
\begin{split}
\left( \frac{d}{2tdt}\right) ^{l+1}P_{2k}\left( \frac{t}{x}\right)
&=(-1)^{k+l+1}\frac{\left(\frac 12\right)_k\left( k+\frac{1}{2}%
\right) _{l+1}}{x^{2l+2}\Gamma (k-l)\left( \frac{1}{2}\right) _{l+1}}
\,_{2}F_{1}\left(
l+1-k,k+l+\frac{3}{2};l+\frac{3}{2};\frac{t^{2}}{x^{2}}\right)\\
&=(-1)^{k+l+1}\frac{\Gamma\left(k+l+\frac 32\right)}{x^{2l+2}\Gamma
(k-l)\Gamma\left(l+\frac 32\right)} \,_{2}F_{1}\left(
l+1-k,k+l+\frac{3}{2};l+\frac{3}{2};\frac{t^{2}}{x^{2}}\right)
\end{split}
\end{equation*}%
when $l<k$. Thus,
\begin{equation}
\begin{split}
K(x,t)&=\frac{\sqrt{\pi }t^{l+1}}{x^{2l+3}\Gamma \left( l+\frac{3}{2%
}\right) }\sum\limits_{k=l+1}^{\infty }(-1)^{k}\beta _{k}(x)\frac{\Gamma \left(k+l+\frac 32\right)}{\Gamma
(k-l)\Gamma\left(l+\frac 32\right)}\,_{2}F_{1}\left( l+1-k,k+l+\frac{3}{2}%
;l+\frac{3}{2};\frac{t^{2}}{x^{2}}\right). \\
&=\frac{\sqrt{\pi }t^{l+1}}{x^{2l+3}\Gamma \left( l+\frac{3}{2}%
\right) }\sum\limits_{m=0}^{\infty }(-1)^{m+l+1}\beta _{m+l+1}(x)\\
&\quad \times \frac{\Gamma\left(m+2l+\frac 52\right)}{\Gamma (m+1)\Gamma\left( l+\frac{3}{2}\right)}
\,_{2}F_{1}\left( -m,m+2l+\frac{5}{2};l+\frac{3}{2};\frac{t^{2}}{x^{2}}%
\right).  \label{Last}
\end{split}
\end{equation}%
The series for $\,_{2}F_{1}(a,b;c;z)$ terminates if either $a$ or $b$ is a
nonpositive integer, in which case the function reduces to a polynomial. In
particular, according to \cite[formula 15.4.6]{AbramowitzStegunSpF},
\begin{equation}
\,_{2}F_{1}(-m,m+\alpha +1+\beta ;\alpha +1;x)= \frac{m!}{(\alpha +1)_{m}}
P_{m}^{(\alpha ,\beta )}(1-2x).  \label{Yakoby}
\end{equation}
Substitution of \eqref{Yakoby} into \eqref{Last} gives us \eqref{K(x,t)}.

Using the asymptotic formula \cite[(6.1.40)]{AbramowitzStegunSpF} one can check
that
\begin{equation}\label{EstimateGammas}
    \log\frac{\Gamma(m+2l+5/2)}{\Gamma(m+l+3/2)} = (l+1) \log (m+l+1) +
    O\left(\frac 1{m+l+1}\right),\qquad m\to\infty.
\end{equation}
Theorem 7.32.2 from \cite{Szego1959} states that
\begin{equation}\label{EstimatePeps}
    \left|P_m^{(l+1/2,l+1)}(z)\right|\le \frac {C_\varepsilon}{\sqrt m}
\end{equation}
uniformly for $z\in [-1+\varepsilon, 1-\varepsilon]$. And it was shown in
\cite[(4.15)]{KTC2017} that for $q\in C^1[0,b]$
\begin{equation}\label{EstimateBeta}
    |\beta_{m+l+1}(x)|\le \frac{c x^{l+3}}{(m+l)^{l+2}},\qquad m\ge 2.
\end{equation}
Combining \eqref{EstimateGammas}, \eqref{EstimatePeps} and \eqref{EstimateBeta}
we obtain the absolute convergence of the series \eqref{K(x,t)} for $x\in(0,b]$
and $t\in (0,x)$, uniform with respect to $t$ on any $[\varepsilon_1,
x-\varepsilon_1]$.

For the whole segment $[0,x]$ note that the Jacobi polynomials satisfy
\begin{equation}\label{EspimatePglob}
    |P_m^{(l+1/2, l+1)}(z)|\le C m^{l+1},\qquad z\in [-1,1],
\end{equation}
see \cite[Theorem 7.32.4]{Szego1959}, while for $q\in C^{2l+5}[0,b]$ the
coefficients $\beta_{m+l+1}$ satisfy \cite[(4.15)]{KTC2017}
\begin{equation}\label{EstimateBeta2}
    |\beta_{m+l+1}(x)|\le \frac{c_2 x^{2l+5}}{(m+l)^{2l+4}}.
\end{equation}
Combining \eqref{EstimateGammas}, \eqref{EspimatePglob} with
\eqref{EstimateBeta2} we obtain the uniform convergence with respect to $t$ on
the whole segment $[0,x]$.

The possibility of termwise differentiation of the series \eqref{Series1}
follows directly from the presented results. Namely, the expressions
$\left(\frac d{2t\,dt}\right)^j R(x,t)$, $j\le l+1$ lead to series similar to
\eqref{K(x,t)} but having $j$ instead of $l+1$, uniformly convergent with
respect to $t$ on any segment $[\varepsilon, x-\varepsilon]$. \hfill\qed

\begin{corollary}\label{Cor3}
Let $q\in C^{2l+5}[0,b]$. Using the formula
\begin{equation*}
P_{n}^{{(\alpha ,\beta )}}(-1)=(-1)^{n}{\binom{n+\beta }{n},}
\end{equation*}%
from \eqref{K(x,t)} and taking into account \eqref{Goursat} we find the
relation
\begin{equation}
\frac 12\int\limits_0^x q(s)\,ds=\frac{(-1)^{l+1}}{x^{l+2}}\frac{\sqrt{\pi }}{\Gamma \left( l+\frac{3}{%
2}\right) }\sum\limits_{m=0}^{\infty }\beta _{m+l+1}(x)\frac{\Gamma \left(
m+2l+\frac{5}{2}\right) }{\Gamma \left( m+l+\frac{3}{2}\right) }{\binom{m+l+1%
}{m}}.  \label{K(x,x)}
\end{equation}
\end{corollary}

The representation \eqref{K(x,t)} may be substituted termwise in
\eqref{transmutation operator} under less restrictive convergence assumption
than the uniform convergence. Namely, $L_1[0,x]$ convergence with respect to
$t$ is sufficient to apply the transmutation operator \eqref{transmutation
operator} with the integral kernel given by \eqref{K(x,t)} to a bounded
function.  Let us rewrite the formula \eqref{K(x,t)} as
\begin{equation}
K(x,t)=\frac{(-1)^{l+1}\sqrt{\pi}}{x^{2l+3}\Gamma \left( l+\frac{3}{2}\right) }
\sum\limits_{m=0}^{\infty }\frac{(-1)^{m}\Gamma \left( m+2l+\frac{5}{2}%
\right) }{\Gamma \left( m+l+\frac{3}{2}\right) } \beta
_{m+l+1}(x)t^{l+1}P_{m}^{\left( l+\frac{1}{2},l+1\right) }\left(
1-2\frac{t^{2}}{x^{2}}\right). \label{K(x,t)L1}
\end{equation}
The following proposition states the $L_1$ convergence of the series
\eqref{K(x,t)L1} under slightly relaxed requirement on the smoothness of the
potential $q$.

\begin{proposition}\label{Prop L1}
Let $l\ge 1$. Under the condition $q\in C^{2l+1}[0,b]$ the series
\eqref{K(x,t)L1} converges for every fixed $x$ in $L_1[0,x]$ norm.
\end{proposition}

\noindent\textbf{Proof.} Consider
\begin{equation}\label{pm}
    \begin{split}
    p_m:&=\int\limits_0^x t^{l+1}\left|P_m^{(l+1/2,
    l+1)}\left(1-2\frac{t^2}{x^2}\right)\right|\,dt =
    \frac{x^{l+2}}{2^{l/2+2}}\int\limits_{-1}^1 (1-z)^{l/2} \left|P_m^{(l+1/2,
    l+1)}(z)\right|\,dz \\
    &\le \frac{x^{l+2}}{2^{l/2+2}}\int\limits_{0}^1 (1-z)^{l/2}
    \left|P_m^{(l+1/2, l+1)}(z)\right|\,dz + \frac{x^{l+2}}{4}\int\limits_{0}^1
    \left|P_m^{(l+1, l+1/2)}(z)\right|\,dz.
    \end{split}
\end{equation}
Theorem 7.34 from \cite{Szego1959} states that
\[
\int\limits_0^1 (1-x)^\mu \left| P_n^{(\alpha,\beta)}(x)\right| \,dx \sim
\begin{cases}
n^{\alpha-2\mu-2}, &\text{if } 2\mu<\alpha-\frac 32,\\
n^{-1/2},&\text{if } 2\mu>\alpha- 3/2,
\end{cases}
\qquad n\to \infty,
\]
whenever $\alpha$, $\beta$, $\mu$ are real numbers greater than $-1$. Using
this result to estimate both integrals in \eqref{pm} we obtain that for some
constant $C$ and all $m\ge 1$
\[
p_m \le Cx^{l+2} m^{l-1}.
\]
Hence $L_1[0,x]$ norms of the functions
$\frac{\Gamma(m+2l+5/2)}{\Gamma(m+l+3/2)}P_m^{(l+1/2,l+1)}(1-2\frac{t^2}{x^2})$
grow at most as $m^{2l}$, $m\to \infty$, and the estimate
\[
|\beta_{m+l+1}(x)|\le \frac{c_3x^{2l+3}}{(m+l)^{2l+2}},
\]
valid (see \cite[(4.150]{KTC2017}) under the condition $q\in C^{2l+1}[0,b]$, is
sufficient to assure the convergence of the series \eqref{K(x,t)} in the
$L_1[0,x]$ norm. \hfill\qed

\begin{remark}
The smoothness requirements $q\in C^{2l+5}[0,b]$  in Theorem \ref{Thm2} and Corollary \ref{Cor3} and $q\in C^{2l+1}[0,b]$ in Proposition \ref{Prop L1} may be excessive. However the minimal requirement $q\in C^1[0,b]$ may be insufficient in general neither for the representation \eqref{K(x,x)} to converge, nor for \eqref{K(x,t)} to converge in $L_1[0,x]$. Indeed, one can easily check that the factor $\frac{\Gamma \left(
m+2l+5/2\right) }{\Gamma \left( m+l+3/2\right) }{\binom{m+l+1%
}{m}}$ grows as $(m+l+1)^{2l+2}$, $m\to\infty$ requiring the coefficients
$\beta_{m}$ to decay faster than $m^{-2l-2}$ in order to fulfill at least the
necessary convergence condition (terms of a series goes to zero as
$m\to\infty$). Similarly, by slightly changing the reasoning in the proof of
Proposition \ref{Prop L1} one can see that the numbers $p_m$ grow as $m^{l-1}$,
requiring the coefficients $\beta_{m}$ to decay faster than $m^{-2l}$ for the
$L_1$ convergence of the series \eqref{K(x,t)}.

Numerical experiments similar to those from \cite[Section 9.1]{KTC2017} suggest
that this may not happen. For the potential
\[
q_2(x)=\begin{cases}
1, & x\in [0,\pi/2],\\
1+(x-\pi/2)^2, & x\in [\pi/2,\pi],
\end{cases}
\]
having its second derivative bounded on $[0,\pi]$, and for $l=3$, the observed
decay rate of the coefficients $\beta_m(x)$ was $Cm^{-7.5}$, $m\to\infty$,
insufficient for the convergence of the series \eqref{K(x,x)}. While for $l=5$
the observed decay rate of the coefficients $\beta_m(x)$ was $Cm^{-9.5}$,
$m\to\infty$, insufficient for the $L_1$ convergence of the series
\eqref{K(x,t)}.
\end{remark}

\begin{corollary} The integral of the power function multiplied by the
kernel $K(x,t)$ has the form
\begin{equation}
\begin{split}
\int\limits_{0}^{x}t^{\alpha }K(x,t)dt&=
\frac{x^{\alpha -l-1}}{2\left( \frac{%
\alpha +l}{2}+1\right) }\frac{\sqrt{\pi }}{\Gamma ^{2}\left( l+\frac{3}{2}%
\right) }\sum\limits_{m=0}^{\infty }(-1)^{m+l+1}\beta _{m+l+1}(x)
\frac{%
\Gamma \left( m+2l+\frac{5}{2}\right) }{m!}\\
&\quad \times \,_{3}F_{2}\left({-}m,2l{+}m{+}\frac{5}{2},\frac{\alpha{+}l}{2}{+}1;l{+}\frac{3%
}{2},\frac{\alpha{+}l}{2}{+}2;1\right),\qquad \operatorname{Re}\alpha>-l-2.
\label{Power}
\end{split}
\end{equation}
\end{corollary}
\noindent\textbf{Proof.}
 Due to \eqref{Last}, we have
\begin{equation}
\begin{split}
\int\limits_{0}^{x}t^{\alpha }K(x,t)dt&=\frac{1}{x^{2l+3}}\frac{\sqrt{\pi }}{%
\Gamma \left( l+\frac{3}{2}\right) } \sum\limits_{m=0}^{\infty
}(-1)^{m+l+1}\beta _{m+l+1}(x)\frac{\Gamma
\left( m+2l+\frac{5}{2}\right) }{\Gamma \left( m+l+\frac{3}{2}\right) }\\
& \quad \times\int\limits_{0}^{x}t^{\alpha +l+1}P_{m}^{\left(
l+\frac{1}{2},l+1\right) }\left( 1-2\frac{t^{2}}{x^{2}}\right) dt.
\label{IntLem}
\end{split}
\end{equation}%
Consider the integral
\begin{equation*}
\int\limits_{0}^{x}t^{\alpha +l+1}P_{m}^{\left( l+\frac{1}{2},l+1\right)
}\left( 1-2\frac{t^{2}}{x^{2}}\right) dt=2^{-2-\frac{\alpha +l}{2}}x^{\alpha +l+2}\int\limits_{-1}^{1}(1-z)^{\frac{%
\alpha +l}{2}}P_{m}^{\left( l+\frac{1}{2},l+1\right) }(z)dz.
\end{equation*}%
Here we use the formula 16.4.3 from  \cite{Bateman} \footnote{Unfortunately,
the formulas 7.391.2 from \cite{Gradshteyn} and 2.22.2.8 from \cite{Prudnikov}
for the same integral, as well as the formula 16.4.3 in the English edition of
\cite{Bateman} contain mistakes.}
\begin{equation*}
\begin{split}
\int\limits_{-1}^{1}(1-x)^{\rho }(1+x)^{\sigma }P_{n}^{(\alpha ,\beta )}(x)dx&=
\frac{2^{\rho +\sigma +1}\Gamma (\rho +1)\Gamma (\sigma +1)}{\Gamma
(\rho +\sigma +2)}\,\frac{\Gamma (n+\alpha +1)}{n!\Gamma (\alpha +1)}\\
&\quad \times \,_{3}F_{2}(-n,\alpha +\beta +n+1,\rho +1;\alpha +1,\rho +\sigma
+2;1),\\
&\mathrm{Re}\,\rho >-1,\qquad\mathrm{Re}\,\sigma >-1.
\end{split}
\end{equation*}%
In particular, we have
\begin{equation}
\begin{split}
\int\limits_{-1}^{1}(1-z)^{\frac{\alpha +l}{2}}P_{m}^{\left( l+\frac{1}{2}%
,l+1\right) }(z)dz&=
\frac{2^{\frac{\alpha +l}{2}+1}}{\left( \frac{\alpha +l}{2}+1\right) }\frac{%
\Gamma \left( m+l+\frac{3}{2}\right) }{m!\Gamma \left( l+\frac{3}{2}\right) }\\
&\quad \times\,_{3}F_{2}\left({-}m,2l{+}m{+}\frac{5}{2},\frac{\alpha{+}l}{2}{+}1;l{+}\frac{3}{2},%
\frac{\alpha{+}l}{2}{+}2;1\right) .  \label{PrInt}
\end{split}
\end{equation}%
Substitution of \eqref{PrInt} into \eqref{IntLem} gives us \eqref{Power}.
\hfill\qed

\begin{remark} In \cite{Castillo1} the images of $x^{2k+l+1}$ under the
action of the operator $T$ were obtained for any $l\geq -1/2$ and $%
k=0,1,2,\ldots $. In the case of integer $l$ they can be derived from %
\eqref{Power}. For example, one can see using \eqref{Power} together with
the equality \cite[formula 9.1]{KTC2017} $\sum\limits_{k=0}^{%
\infty }\beta _{k}(x)=0$ that for $\alpha =1$ and $l=0$ the following relation
is valid
\begin{equation*}
\int\limits_{0}^{x}tK(x,t)dt=\beta _{0}(x).
\end{equation*}
\end{remark}

\section{Representation of the regular solution for
$l=1,2,\ldots$}\label{Sect4}
Substituting \eqref{K(x,t)L1} into \eqref{transmutation operator} we obtain
that the regular solution of \eqref{Eq01} has the form
\begin{equation*}
\begin{split}
    u(\omega, x) &= \sqrt x J_{l+1/2}(\omega x) +
    \frac{\sqrt\pi}{x^{2l+3}\Gamma(l+3/2)} \sum_{m=0}^\infty
    \frac{(-1)^{m+l+1}\Gamma(m+2l+5/2)}{\Gamma(m+l+3/2)}\beta_{m+l+1}(x) \\
    &\quad \times \int\limits_0^x t^{l+3/2}J_{l+1/2}(\omega t)P_m^{(l+1/2,
    l+1)}\left( 1-2\frac{t^2}{x^2}\right)\,dt.
\end{split}
\end{equation*}
Consider the integrals
\[
\mathcal{I}_{k,m}(\omega, x) := \int\limits_0^x t^{k+3/2}J_{k+1/2}(\omega
t)P_m^{(k+1/2, k+1)}\left( 1-2\frac{t^2}{x^2}\right)\,dt.
\]
For $m=0$ the formula 1.8.1.3 from \cite{Prudnikov} gives
\begin{equation}\label{Int0}
\mathcal{I}_{k,0}(\omega, x) = \frac{x^{k+3/2}}{\omega}J_{k+3/2}(\omega x).
\end{equation}
Let $m>0$. Integrating by parts, taking into account \eqref{Int0} and noting
that $\frac d{dx}P_n^{(\alpha,\beta)}(x) = \frac
12(n+\alpha+\beta+1)P_{n-1}^{(\alpha+1,\beta+1)}(x)$ we obtain that
\begin{equation*}
    \begin{split}
    \mathcal{I}_{k,m}(\omega, x) & = \left.
    \frac{t^{k+3/2}}{\omega}J_{k+3/2}(\omega t) P_m^{(k+1/2,
    k+1)}\left(1-2\frac{t^2}{x^2}\right)\right|_{t=0}^x \\
    &\quad - \int\limits_0^x \frac{t^{k+3/2}}{\omega} J_{k+3/2}(\omega t)\cdot
    \frac 12
    (m+2k+5/2)P_{m-1}^{(k+3/2,k+2)}\left(1-2\frac{t^2}{x^2}\right)\cdot\left(-\frac{4t}{x^2}\right)\,dt\\
    &= (-1)^m \binom{m+k+1}{m}\mathcal{I}_{k,0}(\omega, x) +
    \frac{2m+4k+5}{\omega x^2} \mathcal{I}_{k+1, m-1}(\omega, x).
    \end{split}
\end{equation*}
Hence we obtain the following result.
\begin{theorem}\label{Thm repres}
Let $q\in C^{2l+1}[0,b]$. Then the regular solution $u(\omega, x)$ of equation
\eqref{Eq01} satisfying the asymptotics $u(\omega, x)\sim \frac{(\omega
x)^{l+1}}{2^{l+1/2}\Gamma(l+3/2)}$ when $x\to 0$ has the form
\begin{equation}\label{RegSolRep}
    u(\omega, x) = \sqrt{\omega x} J_{l+1/2}(\omega x) +
    \frac{\sqrt{\pi\omega}}{x^{2l+3}\Gamma(l+3/2)} \sum_{m=0}^\infty
    \frac{(-1)^{m+l+1}\Gamma(m+2l+5/2)}{\Gamma(m+l+3/2)}\beta_{m+l+1}(x)\mathcal{I}_{l,m}(\omega,
    x),
\end{equation}
where the coefficients $\beta_k$ are those from \eqref{Series1} and the
functions $\mathcal{I}_{l,m}$ are given by the following recurrent relations
\begin{align}
    \mathcal{I}_{k,0}(\omega, x) &= \frac{x^{k+3/2}}{\omega}J_{k+3/2}(\omega
    x),\qquad k=l,l+1,\ldots,\\
    \mathcal{I}_{k,m}(\omega, x) &= (-1)^m
    \binom{k+m+1}{m}\mathcal{I}_{k,0}(\omega, x) + \frac{2m+4k+5}{\omega x^2}
    \mathcal{I}_{k+1, m-1}(\omega, x),\quad m\in\mathbb{N}.
\end{align}
\end{theorem}

Since the representation \eqref{RegSolRep} is obtained using the transmutation
operator whose kernel is $\omega$-independent, partial sums of the series
\eqref{RegSolRep} satisfy the following uniform approximation property. Let
\begin{equation}\label{K approx}
    K_N(x,t)=\frac{\sqrt{\pi }t^{l+1}}{x^{2l+3}\Gamma \left( l+\frac{3}{2}\right) }
\sum\limits_{m=0}^{N}\frac{(-1)^{m+l+1}\Gamma \left( m+2l+\frac{5}{2}%
\right) }{\Gamma \left( m+l+\frac{3}{2}\right) } \beta _{m+l+1}(x)P_{m}^{\left(
l+\frac{1}{2},l+1\right) }\left( 1-2\frac{t^{2}}{x^{2}}\right),
\end{equation}
and
\begin{equation}\label{u approx}
\begin{split}
    u_N(\omega, x) &= \sqrt{\omega x} J_{l+1/2}(\omega x) \\
    &\quad + \frac{\sqrt{\pi\omega}}{x^{2l+3}\Gamma(l+3/2)} \sum_{m=0}^N
    \frac{(-1)^{m+l+1}\Gamma(m+2l+5/2)}{\Gamma(m+l+3/2)}\beta_{m+l+1}(x)\mathcal{I}_{l,m}(\omega,
    x).
\end{split}
\end{equation}
Due to the $L_1$ convergence of the series \eqref{K(x,t)} with respect to $t$,
for each $x\in (0,b]$ there exists
\begin{equation}\label{Error K Kn}
    \varepsilon_N(x) = \|K(x,\cdot) - K_N(x,\cdot)\|_{L_1[0,x]}
\end{equation}
and $\varepsilon_N(x)\to 0$ as $N\to \infty$.

Due to the asymptotic expansion of the Bessel function for large arguments
\cite[9.2.1]{AbramowitzStegunSpF} there exists a constant $c_l$ such that
\begin{equation}\label{Estimate Bessel}
    |\sqrt z J_{l+1/2}(z)| \le c_l, \qquad z\in\mathbb{R}.
\end{equation}

\begin{theorem}[(Uniform approximation property)]\label{Thm uniform}
Under the conditions of Proposition \ref{Prop L1} the following estimate holds
\begin{equation}\label{u estimate}
    |u(\omega, x) - u_N(\omega, x)|\le c_l \varepsilon_N(x)
\end{equation}
for any $\omega\in \mathbb{R}$.
\end{theorem}

\noindent\textbf{Proof.} Since the functions $u(\omega, x)$ and $u_N(\omega,
x)$ are the images of the same function $\sqrt{\omega x} J_{l+1/2}(\omega x)$
under the action of the integral operators of the form \eqref{transmutation
operator}, one with the integral kernel $K$ and second with the integral kernel
$K_N$, we have
\begin{equation*}
    |u(\omega, x) - u_N(\omega, x)|\le \int\limits_0^x |K(x,t) - K_N(x,t)|
    \cdot |\sqrt{\omega t}J_{l+1/2}(\omega t)|\,dt \le c_l \varepsilon_N(x),
    \end{equation*}
where we used \eqref{Error K Kn} and \eqref{Estimate Bessel}. \hfill\qed

\begin{remark}\label{Rem uniform}
Another uniform approximation property was proved in \cite{KTC2017} for the
regular solution $\tilde u(\omega, x)$ of equation \eqref{Eq01} satisfying the
asymptotics $\tilde u(\omega, x)\sim x^{l+1}$, $x\to 0$. Namely, an estimate of
the form
\[
|\tilde u(\omega, x) - \tilde u_N(\omega, x)|\le \sqrt x
\tilde\varepsilon_N(x)
\]
independent on $\omega\in\mathbb{R}$ was proven. The estimate provided by
Theorem \ref{Thm uniform} is better for large values of $\omega$ due to the
following. Since
\[
\tilde u(\omega, x) = \frac{2^{l+1/2}\Gamma(l+3/2)}{\omega^{l+1}} u(\omega,
x),
\]
and for each fixed $x$ the function $u(\omega, x)$ remains bounded as
$\omega\to\infty$, the function $\tilde u(\omega, x)$ decays at least as
$\omega^{-l-1}$ as $\omega\to\infty$, meaning that the uniform error estimate
is useful only in some neighbourhood of $\omega =0$. Meanwhile the estimate
\eqref{u estimate} remains useful even for large values of $\omega$.
\end{remark}

\section{The case of a noninteger $l$}\label{Sect5}

\begin{theorem} Let $q\in C^1[0,b]$. For $l\geq -1/2$ the following formula for
the kernel $K$
is valid
\begin{equation}
K(x,t)=  \frac{\sqrt{\pi}}{\Gamma(l+3/2)}\frac{t^{l+1}}{x(x^2-t^2)^{l+1}}\sum_{k=0}^\infty \frac{(-1)^kk!}{\Gamma(k-l)}\beta_k(x) P_k^{(l+1/2,-l-1)}\left(1-\frac{2t^2}{x^2}\right), \label{K for noninteger}
\end{equation}%
where $P_k^{(l+1/2, -l-1)}$ are polynomials given by the same formulas as the
classical Jacobi polynomials \cite[4.22]{Szego1959},
\cite[22.3]{AbramowitzStegunSpF}.

The series in \eqref{K for noninteger} converges absolutely for any $x\in
(0,b]$ and $t\in (0,x]$, uniformly with respect to $t$ on any segment $t\in
[\varepsilon, x]$. Under the additional assumption that $q\in C^{2[l]+5}[0,b]$,
where $[l]$ denotes the largest integer not exceeding $l$, the convergence is
uniform with respect to $t$ on $[0,x]$.
\end{theorem}

\noindent\textbf{Proof.} Let $l=[l]+\{l\}$, $p:=[l]$, $p=-1,0,1,\ldots$ and
$\lambda :=\{l\}$, $0\leq \lambda <1$.

Taking in formula \eqref{Int4} $n=p+2$ we obtain
\begin{equation}
K(x,t)=\frac{\sqrt{\pi }}{\Gamma \left( l+\frac{3}{2}\right) }\frac{t^{l+1}}{%
\Gamma (1-\lambda )}\left( -\frac{d}{2tdt}\right) ^{p+2}\int\limits_{t}^{x}(s^{2}-t^{2})^{-\lambda }sR(x,s)ds.  \label{NI}
\end{equation}%
Substituting the series expansion
$$R(x,t)=\sum\limits_{k=0}^{\infty }\frac{%
\beta _{k}(x)}{x}P_{2k}\left( \frac{t}{x}\right)$$
 into \eqref{NI} we get (the possibility to differentiate termwise follows
 similarly to the proof of Theorem \ref{Thm2})
\begin{equation}\label{k for noninteger2}
K(x,t)=\frac{2\sqrt{\pi }}{\Gamma \left( l+\frac{3}{2}\right) }\frac{t^{l+1}%
}{\Gamma (1-\lambda )}\sum\limits_{k=0}^{\infty }\frac{\beta _{k}(x)}{x}\left( -\frac{d}{2tdt}\right)
^{p+2}\int\limits_{t}^{x}(s^{2}-t^{2})^{-\lambda }sP_{2k}\left( \frac{s}{x}%
\right) ds.
\end{equation}%
Consider the integral
\[
I_k:= \int\limits_t^x (s^2-t^2)^{-\lambda} s P_{2k}\left(\frac sx\right)\,ds.
\]
Using formula 2.17.2.9 from \cite{Prudnikov} and noting that $B(1-\lambda,
\lambda-1)=0$ we obtain that
\[
I_k=\frac{(-1)^k (\lambda-1/2)_k
x^{2-2\lambda}}{2(1-\lambda)_{k+1}}\,_2F_1\left(\lambda-1-k,\lambda-\frac
12+k;\lambda -\frac 12; \frac{t^2}{x^2}\right).
\]
Proceeding as in the proof of Theorem \ref{Thm2} we get
\[
\begin{split}
\left(-\frac{d}{2t\,dt}\right)^{p+2} I_k &= \frac{(-1)^{k+p+2} (\lambda-1/2)_k
x^{2-2\lambda}}{2(1-\lambda)_{k+1}}
\frac{(\lambda-1-k)_{p+2} (\lambda-1/2+k)_{p+2}}{(\lambda-1/2)_{p+2}
x^{2(p+2)}}\\
&\quad \times\,_2F_1\left(\lambda -k+p+1, \lambda+k+p+\frac 32; \lambda+p+\frac
32; \frac{t^2}{x^2}\right).
\end{split}
\]
Noting that $\lambda+p=l$ and using $(x)_n=\frac{\Gamma(x+n)}{\Gamma(x)}$ we
obtain
\[
\begin{split}
\left(-\frac{d}{2t\,dt}\right)^{p+2} I_k &=
\frac{(-1)^{k+p+2}\Gamma(l+1-k)\Gamma(l+3/2+k)\Gamma(1-\lambda)}{2\Gamma(2+k-\lambda)\Gamma(\lambda-1-k)\Gamma(l+3/2)x^{2l+2}}\\
&\quad\times\,_2F_1\left(l+1-k,l+\frac 32+k; l+\frac 32;
\frac{t^2}{x^2}\right).
\end{split}
\]
Using the reflection formula $\Gamma(z)\Gamma(1-z)=\frac{\pi}{\sin\pi z}$ we
see that
\[
\begin{split}
\frac{\Gamma(l+1-k)}{\Gamma(2+k-\lambda)\Gamma(\lambda-1-k)} &=
\frac{\sin\pi(\lambda-1-k)}{\sin\pi(l+1-k)\Gamma(k-l)} \\
&= \frac{(-1)^{1+k+p}\sin\pi(\lambda + p)}{(-1)^{1-k}\sin\pi
l\,\Gamma(k-l)}=\frac{(-1)^p}{\Gamma(k-l)},
\end{split}
\]
and hence from the formula 15.3.3 from \cite{AbramowitzStegunSpF} we get
\begin{equation}\label{IntegralsNonInt}
\begin{split}
\left(-\frac{d}{2t\,dt}\right)^{p+2} I_k &=
\frac{(-1)^{k}\Gamma(l+3/2+k)\Gamma(1-\lambda)}{2
\Gamma(k-l)\Gamma(l+3/2)x^{2l+2}
\left(1-\frac{t^2}{x^2}\right)^{l+1}}\,_2F_1\left(k+\frac 12,-k; l+\frac 32;
\frac{t^2}{x^2}\right)\\
&=
\frac{(-1)^k\Gamma(l+3/2+k)\Gamma(1-\lambda)}{2\Gamma(k-l)\Gamma(l+3/2)(x^2-t^2)^{l+1}}\,_2F_1\left(-k,k+\frac
12; l+\frac 32; \frac{t^2}{x^2}\right)\\
&= \frac{(-1)^k\Gamma(1-\lambda)k!}{2\Gamma(k-l)
(x^2-t^2)^{l+1}}P_k^{(l+1/2,-l-1)}\left(1-\frac{2t^2}{x^2}\right),
\end{split}
\end{equation}
where $P_n^{(\alpha,\beta)}$ stands for the Jacobi polynomials (see
\cite[15.4.6]{AbramowitzStegunSpF}), however note that due to the second
parameter equal to $-l-1$ the polynomials $P_k^{(l+1/2,-l-1)}$ are not
classical orthogonal polynomials, even though they are given by the same
formulas and satisfy the same recurrence relations, see  \cite[4.22]{Szego1959}
for additional details.

Combining \eqref{k for noninteger2} with \eqref{IntegralsNonInt} finishes the
proof of \eqref{K for noninteger}. Convergence of the series can be obtained
similarly to the proof of Theorem \ref{Thm2} with the only difference that for
$q\in C^{2[l]+5}[0,b]$ the coefficients $\beta_k$ satisfy $|\beta_k(x)|\le
\frac{cx^{l+2}}{k^{2l+3}}$, see \cite[(4.17)]{KTC2017}. \hfill\qed

\begin{remark}
The representation \eqref{K for noninteger} is also valid for integer values of
$l$ (and in such case can be simplified to those of \eqref{K(x,t)}). Indeed,
for integer $l\ge 0$ the terms having $k\le l$ in \eqref{K for noninteger} are
equal to zero due for the factor $\Gamma(k-l)$. While for $k\ge l+1$ one
obtains using the formula 4.22.2 from \cite{Szego1959} that
\begin{equation}
\begin{split}
P_k^{(l+1/2, -l-1)}(x) &= (-1)^k P_k^{(-l-1,l+1/2)}(-x) \\
&= \frac{(-1)^k}{\binom{k}{l+1}}\binom{k+l+1/2}{l+1}
\left(\frac{-x-1}{2}\right)^{l+1} P_{k-l-1}^{(l+1, l+1/2)}(-x) \\
&= \frac{\Gamma(k+l+3/2) (k-l-1)!}{\Gamma(k+1/2) k!}
\left(\frac{x+1}{2}\right)^{l+1} P_{k-l-1}^{(l+1/2,l+1)}(x).
\end{split}\label{ReduceJacobi}
\end{equation}
Applying \eqref{ReduceJacobi} in \eqref{K for noninteger} one easily arrives at
\eqref{K(x,t)}.
\end{remark}

\section{Numerical illustration}\label{Sect6}

\subsection{Integer $l$: a spectral problem}
The approximate solution \eqref{u approx} can be used for numerical solution of
the Dirichlet spectral problem for equation \eqref{Eq01}, i.e., for finding
those $\omega$ for which there exists a regular solution of \eqref{Eq01}
satisfying \[ u(\omega, b)=0.
\]
The uniform approximation property \eqref{u estimate} leads to a uniform error
bound for both lower and higher index eigenvalues. An algorithm is
straightforward: one computes coefficients $\beta_k$, chooses $N$ as an index
where the values $|\beta_k(b)|$ cease to decay due to machine precision
limitation and looks for zeros of the analytic function $F(\omega):=
u_N(\omega, b)$. We refer the reader to \cite{KTC2017} for implementation
details regarding the computation of $\beta_k$. We want to emphasize that the
presented numerical results are only ``proof of concept'' and are not aimed to
compete with the best existing software packages.

\begin{figure}[htb!]
\centering
\begin{tabular}{ccc}
$l=1$, $N=13$ & & $l=2$, $N=14$ \\
\includegraphics[bb=0 0 180 144, width=2.5in,height=2in]
{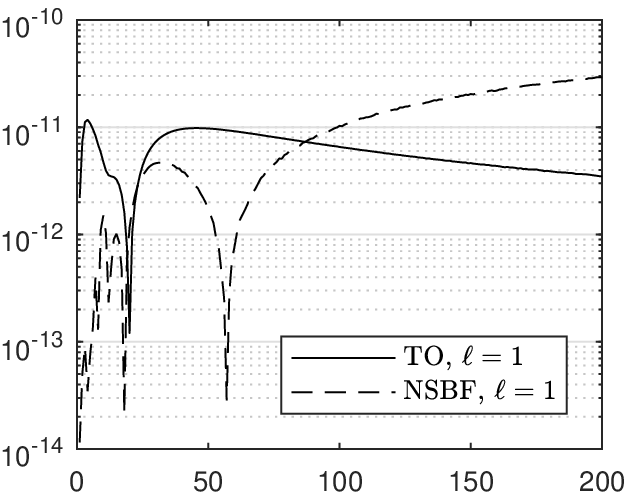} & \ &
\includegraphics[bb=0 0 180 144, width=2.5in,height=2in]
{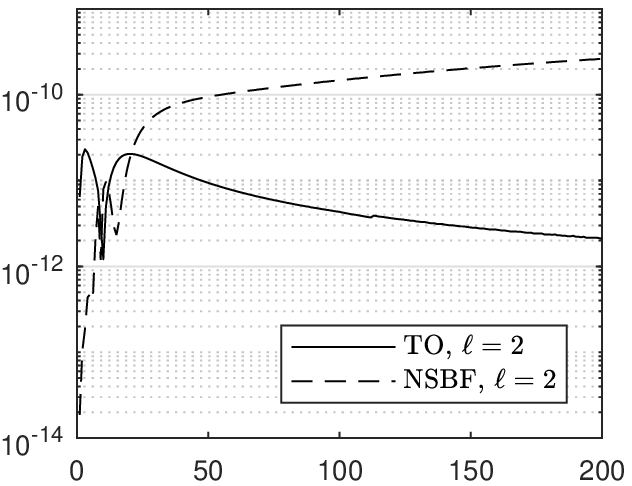}\\
$l=5$, $N=16$ & & $l=10$, $N=19$\\
\includegraphics[bb=0 0 180 144, width=2.5in,height=2in]
{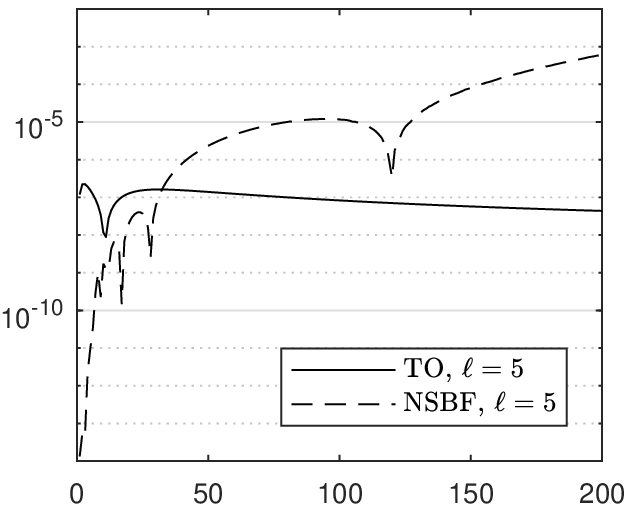} & \ &
\includegraphics[bb=0 0 180 144, width=2.5in,height=2in]
{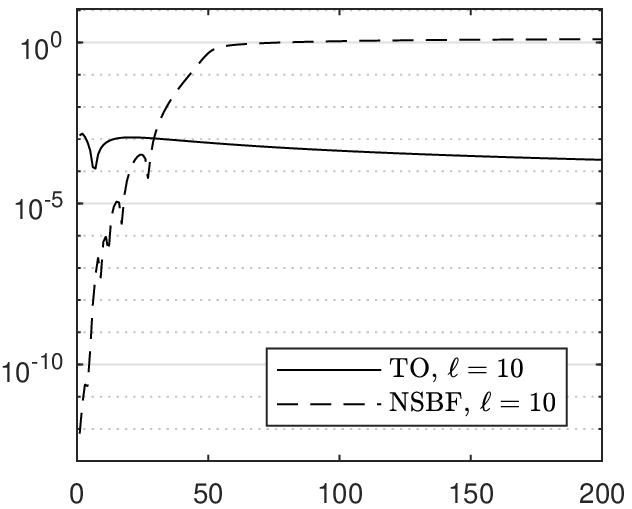}
\end{tabular}
\caption{Absolute errors of the first 200 eigenvalues for the spectral problem
\eqref{Ex1Eq}, \eqref{Ex1BC} for different values of $l$ obtained using the
proposed method (solid lines, marked as TO on the legends) and those produced
by the method from \cite{KTC2017} (dashed lines, marked as NSBF on the
legends). The parameter $N$ over each figure corresponds to the truncation
parameter used for the approximate solution \eqref{u approx}.} \label{Ex1Fig1}
\end{figure}

Consider the following spectral problem
\begin{gather}
-u''+\left(\frac{l(l+1)}{x^2}+x^2\right) u=\omega^2 u, \quad 0\le x\le
\pi,\label{Ex1Eq}\\
u(\omega, \pi)=0.\label{Ex1BC}
\end{gather}
The regular solution of equation \eqref{Ex1Eq} can be written as
\[
u(\omega, x) = x^{l+1}e^{x^2/2}\,_1F_1\left(\frac{\omega^2+2l+3}{4}; l+\frac
32; -x^2\right)
\]
allowing one to compute with any precision arbitrary sets of eigenvalues
using, e.g., Wolfram Mathematica. We compare the results provided by the
proposed algorithm to those of \cite{KTS} and \cite{KTC2017} where other
methods based on the transmutation operators were implemented. The following
values of $l$ were considered: $1$, $2$, $5$ and $10$. All the computations
were performed in machine precision using Matlab 2012. For each value of $l$ we
computed 200 approximate eigenvalues. In Table \ref{Ex1Table1} we show absolute
errors of some eigenvalues for $l=1$. For $l>1$ the analytic approximation
proposed in \cite{KTS} produced considerably worse results, for that reason we
compared the approximate results only with those from \cite{KTC2017}. We
present the results on Figure \ref{Ex1Fig1}.

\begin{table}[htb!]
\centering
\begin{tabular}{c|c|c|c|c}\hline
$n$ & $\omega_{n}$ (Exact)  & $\Delta \omega_n$ \eqref{u approx}& $\Delta
\omega_n$ (\cite{KTC2017}) & $\Delta \omega_n$ (\cite{KTS}) \\\hline
1 &  $2.24366651120741$ & $2.2\cdot 10^{-12}$ & $1.1\cdot 10^{-14}$ & $3.6\cdot
10^{-6}$ \\
2 &  $3.09030600792814$ & $7.4\cdot 10^{-12}$ & $5.1\cdot 10^{-14}$ & $4.1\cdot
10^{-6}$ \\
5 &  $5.78188700721372$ & $1.1\cdot 10^{-11}$ & $7.6\cdot 10^{-14}$ & $2.1\cdot
10^{-5}$ \\
10 & $10.6472529934013$ &  $4.7\cdot 10^{-12}$ & $1.5\cdot 10^{-12}$ &
$9.9\cdot 10^{-6}$\\
20 & $20.5753329357456$ & $1.2\cdot 10^{-13}$  & $1.2\cdot 10^{-12}$ &
$3.1\cdot 10^{-6}$\\
50 & $50.5305689586825$ &  $9.7\cdot 10^{-12}$ & $1.7\cdot 10^{-12}$ &
$3.7\cdot 10^{-6}$\\
100 &$100.515359633269$ & $6.5\cdot 10^{-12}$ & $1.0\cdot 10^{-11}$ & $3.4\cdot
10^{-7}$\\
200 &$200.507698855317$ & $3.5\cdot 10^{-12}$ & $2.9\cdot 10^{-11}$ & $2.4\cdot
10^{-7}$\\
\hline
\end{tabular}
\caption{The eigenvalues for the spectral problem \eqref{Ex1Eq}, \eqref{Ex1BC}
for $l=1$ compared to those produced by the approximation \eqref{u approx} with
$N=13$ and to those reported in \cite{KTC2017} and \cite{KTS}. $\Delta\omega_n$
denotes the absolute error of the computed eigenvalue $\omega_n$.}
\label{Ex1Table1}
\end{table}

The obtained results confirm Remark \ref{Rem uniform}, the proposed method
outperforms the method from \cite{KTC2017} for large index eigenvalues. Uniform
(and even decaying) absolute error of approximate eigenvalues can be
appreciated. The loss of accuracy for large values of $l$ can be easily
explained from the representation \eqref{RegSolRep}. Recall that due to
recurrent formulas used for computation of the coefficients $\beta_m$ (see
\cite[(6.11)]{KTC2017}) the absolute errors of the computed coefficients
$\beta_m$ are slowly growing as $m\to\infty$. The coefficients
$\Gamma(m+2l+5/2)/\Gamma(m+l+3/2)$ also grow as $m\to\infty$. For $l=5$ the
first coefficient is about $4.8\cdot 10^5$, while for $l=10$ the first
coefficient is about $2\cdot 10^{13}$, which explains the loss of accuracy.

\subsection{Non-integer $l$: approximate integral kernel}
We illustrate the representation \eqref{K for noninteger} constructing
numerically the integral kernel $K$. Unfortunately we are not aware of any
single nontrivial potential $q$ for which the integral kernel $K$ is known in a
closed form. In \cite{KTS} an analytic approximation was proposed and revealed
excellent numerical performance for the potential $q(x)=x^2$, $x\in [0,\pi]$
and for $l=-0.5$ or $l=0.5$ (the Goursat data \eqref{Goursat} was satisfied
with an error less than $10^{-12}$ and a large set of eigenvalues was
calculated with absolute errors smaller than $10^{-11}$). It is worth to
mention that for other values of $l$, say $1/3$ or $3/2$, and for other
potentials, the performance of the method from \cite{KTS} was considerably
worse. We consider the same potential and the same values of $l$ to illustrate
the numerical behavior of the representation \eqref{K for noninteger}, using
the approximate kernel $K$ obtained with the method from \cite{KTS} in the role
of an exact one for all the comparisons.

\begin{figure}[htb!]
\centering
\includegraphics[bb=0 0 360 259, width=5in,height=3.6in]{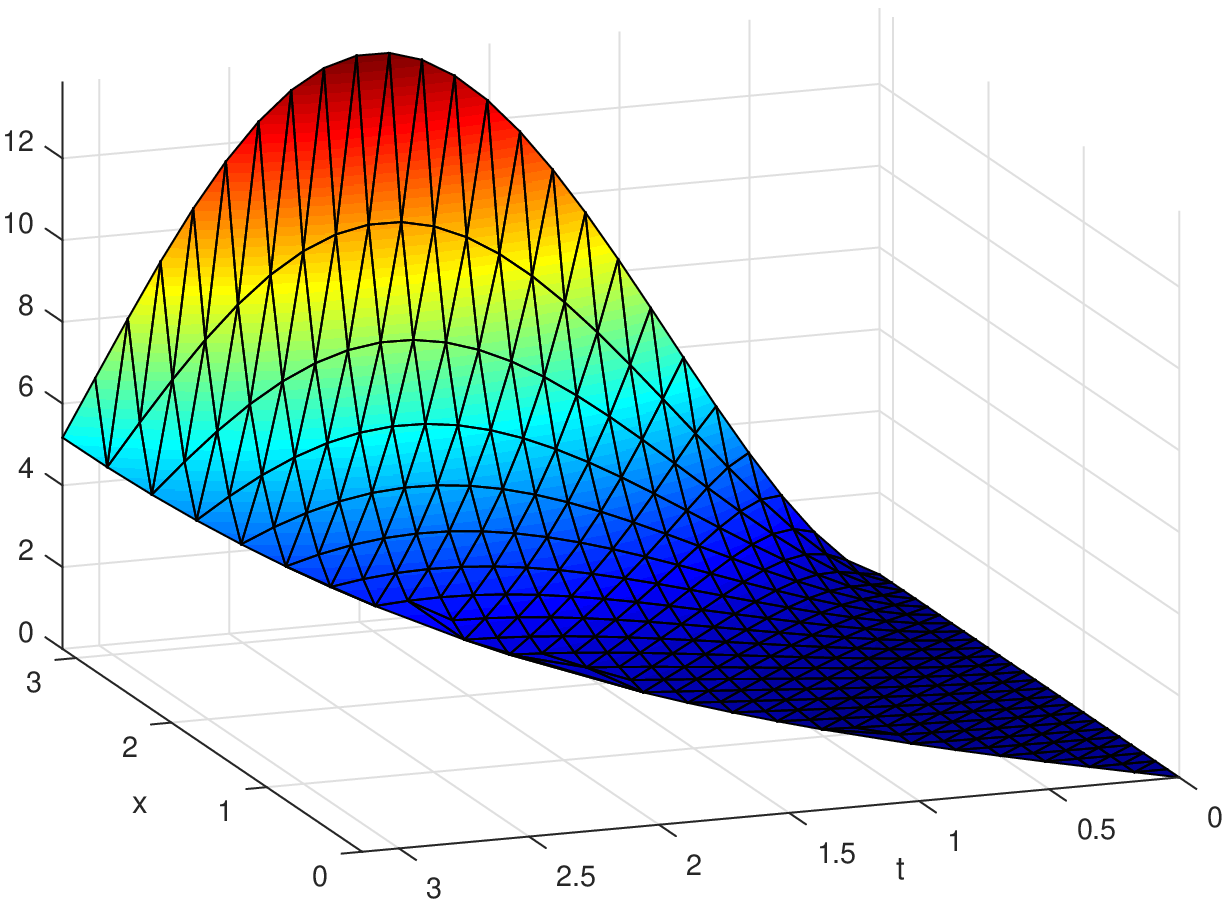}
\caption{The integral kernel $K$ of the transmutation operator for equation
\eqref{Ex1Eq} with $l=0.5$.} \label{Ex2Fig1}
\end{figure}

\begin{figure}[htb!]
\centering
\includegraphics[bb=0 0 360 180, width=5in,height=2.5in]
{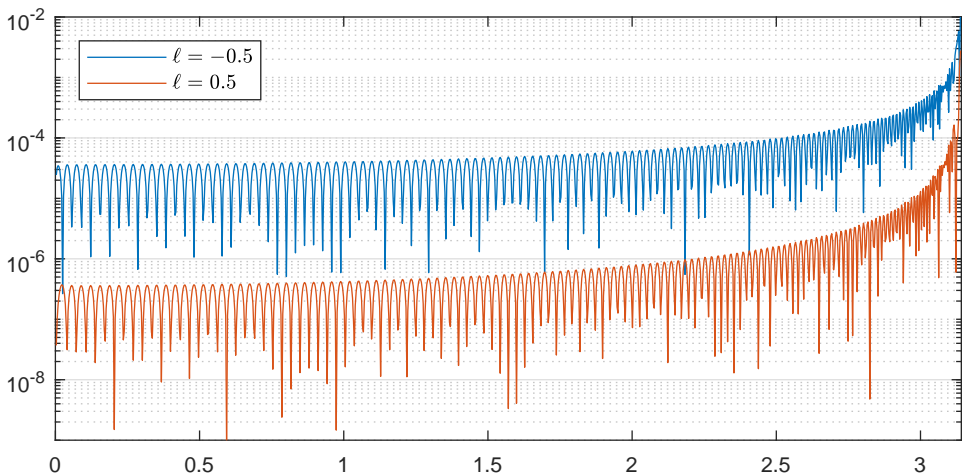} \caption{Absolute error of the approximate integrate kernel for
equation \eqref{Ex1Eq} at $x=\pi$ obtained from the truncated representation
\eqref{K for noninteger} taking coefficients up to $k=150$, i.e., the value
$|K(\pi, t)-K_{150}(\pi, t)|$, $t\in [0,\pi)$. Upper (blue) line corresponds to
$l=-0.5$, lower (red) line corresponds to $l=0.5$.} \label{Ex2Fig2}
\end{figure}

It is known \cite[Proposition 4.5]{KTC2017} that for non-integer values of $l$
the coefficients $\beta_k$ decay as $k^{-2l-3}$ when $k\to\infty$, hence the
series in \eqref{K for noninteger} converges rather slow. For that reason we
computed the coefficients $\beta_k$ for $k\le 150$. The necessary values of the
Jacobi polynomials $P_k^{l+1/2, -l-1}$ were calculated using the recurrent
formula (4.5.1) from \cite{Szego1959}. It is worth to mention that the whole
computation took few seconds. On Figure \ref{Ex2Fig1} we present the kernel $K$
for $l=0.5$. On Figure \ref{Ex2Fig2} we show the absolute error, the value of
the difference $|K(\pi, t)-K_{150}(\pi, t)|$, $t\in [0,\pi)$ for $l=-0.5$ and
$l=0.5$. The growth of the error as  $t\to x$ can be explained by the division
over $(x^2-t^2)^{l+1}$ in \eqref{K for noninteger}. Nevertheless, a remarkable
accuracy can be appreciated. The obtained approximation may be used, in
particular, for solution of spectral problems, we leave the detailed analysis
for a separate study.

\section*{ACKNOWLEDGEMENTS}
Research was supported by CONACYT, Mexico via the project 222478.

\end{document}